\documentstyle[amscd,amssymb,verbatim,12pt]{amsart}
\newcommand{\C}{{\Bbb C}}

\newcommand{\const}{\mbox{\rm const.}}

\newcommand{\diam}{\mbox{\rm diam}}

\newcommand{\Image}{\mbox{\rm Im}}

\newcommand{\Ker}{\mbox{\rm Ker}}

\newcommand{\R}{{\Bbb R}}

\newcommand{\Ric}{\mbox{\rm Ric}}

\newcommand{\vol}{\mbox{\rm vol}}
\newcommand{\Z}{{\Bbb Z}}
\theoremstyle{plain}

\newtheorem{lemma}{Lemma}
\newtheorem{theorem}{Theorem}
\newtheorem{proposition}{Proposition}
\newtheorem{corollary}{Corollary}

\numberwithin{equation}{section}
\renewcommand{\rm}{\normalshape}
\begin{document}
\title{On the Spectrum of a Finite-Volume Negatively-Curved Manifold}
\author{John Lott}
\address{Department of Mathematics\\
University of Michigan\\
Ann Arbor, MI  48109-1109\\
USA}
\email{lott@@math.lsa.umich.edu}
\thanks{Research supported by NSF Grant DMS-9704633}
\date{September 10, 1999}
\maketitle
\begin{abstract}
We show that a noncompact manifold with bounded sectional curvature, 
whose ends are sufficiently Gromov-Hausdorff close to rays,
has a finite dimensional space of square-integrable
harmonic forms.  In the special case of a finite-volume manifold with pinched
negative sectional curvature, we show that the essential spectrum of the
$p$-form Laplacian is the union of the essential spectra of a collection
of ordinary differential operators associated to the ends.  We give
examples of such manifolds with curvature pinched arbitrarily close
to $-1$ and with an infinite number of gaps in the spectrum of the
function Laplacian.
\end{abstract}
\section{Introduction} \label{sect1}

In this paper we consider 
Riemannian manifolds of finite volume and pinched negative sectional curvature.
We give results about the kernel of the differential form
Laplacian and about its essential
spectrum.

Our first result is the finite dimensionality of the space of square-integrable
harmonic forms for a more general class of Riemannian 
manifolds, which can be roughly characterized as those with
bounded sectional curvature and with ends that are sufficiently 
Gromov-Hausdorff close to rays.
Let $M$ be a complete connected $n$-dimensional
Riemannian manifold with a basepoint $m$.
Let $B_r(m)$ denote the distance ball around $m$ and let $S_r(m) \: = \:
\partial B_r(m)$ be the distance sphere around $m$. Put 
\begin{equation}
{\cal D}_r(m) \: = \: \sup_{\Sigma_r} \diam(\Sigma_r),
\end{equation}
where $\Sigma_r$ ranges over the connected components of $S_r(m)$ which
intersect a ray through $m$.

For $p \in \Z \cap [0, n]$,
let $\triangle^M_p$ be the $p$-form Laplacian on $M$. A harmonic $p$-form
on $M$ is an element of $\Ker(\triangle^M_p)$.
Let ${\cal H}^p_{(2)}(M)$ denote the vector space
of square-integrable harmonic $p$-forms on $M$.

\begin{theorem} \label{thm1}
There is a number $\delta \: = \: \delta(n) \: >  \:0$ with
the property that if for some $b > 0$ the sectional curvatures of $M$ are all
bounded in absolute value by $b^2$, and 
\begin{equation} \label{delta}
\limsup_{r \rightarrow \infty} {\cal D}_r(p) \: \le \: \delta \: b^{-1},
\end{equation}
then for all $p \in \Z \cap [0, n]$ the dimension of ${\cal H}^p_{(2)}(M)$
is finite.
\end{theorem}

\begin{corollary} \label{cor1}
Let $M$ be a complete connected $n$-dimensional
Riemannian manifold of finite volume whose
sectional curvatures satisfy $- \: b^2 \: \le \: K \: \le \: - \: a^2$, 
with $0 \: < \: a \: \le \: b$. Then for all $p \in \Z \cap [0, n]$,
the dimension of ${\cal H}^p_{(2)}(M)$
is finite.
\end{corollary}

Corollary \ref{cor1} was previously known to be true when
$p \notin \{ \frac{n - 1}{2}, \frac{n}{2}, \frac{n + 1}{2} \}$ and 
$\frac{b}{a} \: < \: 
\frac{n-1}{2 \min (p, n-p)}$,
and when $p \: = \: \frac{n}{2}$ and $\frac{b}{a}$ satisfies a certain 
inequality for which we refer to \cite{Donnelly-Xavier (1984)}.

The other results in this paper concern manifolds $M$ as in Corollary
\ref{cor1}. Recall that the essential spectrum of $\triangle^M_p$
consists of all numbers in the spectrum of $\triangle^M_p$ 
other than those which are both isolated in the spectrum and have
a finite-dimensional eigenspace.

Label the ends of $M$ by $I \in \{1, \ldots, B\}$.
An end of $M$ has a neighborhood $U_I$ whose closure 
is homeomorphic to 
$[0, \infty) \times N_I$, with $N_I$ an infranilmanifold and the 
parameter $s \in [0, \infty)$ being a Busemann function for the end.
As will be explained, if $U_I$ lies far enough out the end then
the differential forms on each fiber $\{ s \} \times N_I$ decompose into
a finite-dimensional space $E_I(s)$, consisting of ``bounded energy'' forms,
and its orthogonal complement $E_I(s)^\perp$, consisting of ``high energy''
forms.  The vector spaces $\{E_I(s)\}_{s \in [0, \infty)}$ fit together to
form a vector
bundle $E_I$ on $[0, \infty)$. Let $P_0$ be orthogonal projection from
$\bigoplus_{I=1}^B \Omega^*(\overline{U_I})$ to 
$\bigoplus_{I=1}^B \Omega^*([0, \infty); E_I)$.
Put ${\cal A} \: = \: P_0 \: d^M \: P_0$. Consider the operator
${\cal A} \: {\cal A}^* \: + \: 
{\cal A}^* \: {\cal A}$ corresponding to the quadratic form
\begin{equation} \label{Q(F,F)}
Q(\omega) \: = \: \int_0^\infty \left[
|{\cal A} \omega|^2 \: + \: |{\cal A}^* \omega|^2 \right] \: ds,
\end{equation}
where $\omega \in \bigoplus_{I=1}^B \Omega^*([0, \infty); E_I)$ satisfies
the boundary condition that its pullback to $\{0\}$ vanishes.
Then ${\cal A} \: {\cal A}^* \: + \: {\cal A}^* \: {\cal A}$ 
is a second-order ordinary
differential operator. Let
$({\cal A} \: {\cal A}^* \: + \: {\cal A}^* \: {\cal A})_p$ 
denote its restriction
to elements of total degree $p$.

\begin{theorem} \label{thm2}
Suppose that $M$ is as in Corollary \ref{cor1}.
Then for all 
$p \in \Z \cap [0, n]$, the essential spectrum of $\triangle^M_p$ equals
the essential spectrum of $({\cal A} \: {\cal A}^* \: + \: 
{\cal A}^* \: {\cal A})_p$.
\end{theorem}

Theorem \ref{thm2} was previously known in the case when $M$ is a
finite-volume rank-$1$ locally symmetric space \cite{Mueller (1987)}.

As an example of Theorem \ref{thm2}, we consider the case of
the Laplacian on functions. It is well-known that if $M$ is a noncompact
finite-volume hyperbolic manifold then the spectrum of
its function Laplacian is the union of 
$\left[ \frac{(n-1)^2}{4}, \infty \right)$
with a finite subset of $\left[ 0,  \frac{(n-1)^2}{4} \right)$. 
In particular, there is a finite number of gaps in the spectrum.
We show that for manifolds with sectional curvature pinched close to $-1$, the
situation can be very different.

\begin{theorem} \label{thm3} 
For any $\epsilon > 0$, there is a complete connected noncompact 
finite-volume Riemannian manifold
whose sectional curvatures lie in $[-1-\epsilon, -1+\epsilon]$ and whose
function Laplacian has an infinite number of gaps in its spectrum.
The gaps tend toward infinity.
\end{theorem}

Theorems \ref{thm1} and \ref{thm2} continue to hold if one
allows $M$ to be altered within a compact region. The proofs go through
without change.

I am grateful to Lizhen Ji for suggesting this line of research and for
many helpful discussions.  I thank the IHES for its hospitality while
this research was performed.

\section{Proof of Theorem \ref{thm1}} \label{sect2} 

The vector space ${\cal H}^p_{(2)}(M)$ is isomorphic to
the $p$-dimensional (reduced) $L^2$-cohomology of $M$. For 
background on $L^2$-cohomology, we refer to \cite{Lott (1996)}, 
\cite[Section 2]{Lott (1997)} and references
therein.

Suppose that the sectional curvatures of $M$ are all
bounded in absolute value by $b^2$.
From \cite[Theorem 1]{Shen (1994)}, if the number $\delta$ is sufficiently
small and $M$ satisfies (\ref{delta}) then $M$ has
finite topological type, i.e. is homeomorphic to the interior of a 
compact manifold-with-boundary $\overline{M}$. (In fact, for this conclusion
it is enough to just have the lower bound on the sectional curvatures
\cite{Shen (1993)}.) In particular, if $\{N_I\}_{I=1}^B$ are the connected
components of $\partial \overline{M}$ then there is a compact set 
$K \subset M$ such that the closures $\{ \overline{U_I} \}_{I=1}^B$ of
the connected components of
$M - K $ are homeomorphic to $\{ [0, \infty) \times N_I \}_{I=1}^B$.

From \cite[Proposition 5]{Lott (1997)}, the dimension of 
${\cal H}^p_{(2)}(M)$ is finite if and only if the dimension of 
${\cal H}^p_{(2)}(\overline{U_I})$ is finite
for each $I$. Here ${\cal H}^p_{(2)}(\overline{U_I})$ can be defined
either as the $p$-dimensional
(reduced) $L^2$-cohomology of $\overline{U_I}$ or as the
space of square-integrable harmonic $p$-forms on $\overline{U_I}$ satisfying
absolute boundary conditions on $\partial \overline{U_I}$.

From \cite[Theorem 2]{Shen (1994)}, $N_I$ is diffeomorphic to an 
infranilmanifold. The proof of \cite[Theorem 2]{Shen (1994)} uses
the collapsing results of Cheeger, Fukaya 
and Gromov, as given for example in  \cite{Cheeger-Fukaya-Gromov (1992)}.
In particular, it uses Fukaya's
fibration theorem, along with the fact that $U_I$
is Gromov-Hausdorff
close to a ray which passes through it. Strictly speaking, as in the proof of
\cite[Theorem 2]{Shen (1994)}, one may have to shrink $U_I$ a bit
in order to apply the fibration theorem.

In fact, \cite{Cheeger-Fukaya-Gromov (1992)} describes a model metric which
is uniformly $C^0$-close to that of 
$\overline{U_I}$. However, reduced 
$L^2$-cohomology is biLipschitz invariant
(see, for example, \cite[Proposition 1]{Lott (1997)}). Hence it suffices to
compute the (reduced) $L^2$-cohomology of 
$\overline{U_I}$ with the model metric. We now describe the model metric.

The infranilmanifold $N_I$ is $F_I$-covered by a nilmanifold $\Gamma_I
\backslash {\frak N}_I$ where ${\frak N}_I$ is a simply-connected connected
nilpotent Lie group, $\Gamma_I$ is a lattice in ${\frak N}_I$ and
$F_I$ is a finite group of automorphisms of ${\frak N}_I$ which preserve
$\Gamma_I$. 
From \cite[Proposition 4.9]{Cheeger-Fukaya-Gromov (1992)}, the model
metric on $\overline{U_I}$ is that of
a certain Riemannian submersion from $\overline{U_I}$ to
$[0, \infty)$ which is invariant under a local action of ${\frak N}_I$. 
In particular, the flow of the horizontal vector field for the
Riemannian submersion
$\overline{U_I} \rightarrow [0, \infty)$ preserves the affine structures
on the fibers. By integrating the vector field, the model metric can be
written in the form
\begin{equation}
g \: = \: ds^2 \: + \: h(s),
\end{equation}
where for each $s \: \in
\: [0, \infty)$,
$h(s)$ is a smooth 
metric on $N_I$ which comes from an $F_I$-invariant left-invariant
metric on ${\frak N}_I$. Furthermore, if $S(s)$ denotes the second fundamental
form of $\{s \} \times N_I$ then we can assume that $\{S(s)\}_{s \in 
[0, \infty)}$ are uniformly bounded with respect to $\{h(s)\}_{s \in 
[0, \infty)}$.  In what follows we will allow ourselves to reduce the
end by making finite shifts of the interval $[0, \infty)$, without change
of notation. 

There is a canonical flat connection $\nabla^{aff}$ on $TN_I$ coming
from the flat connection on $T{\frak N}_I$ for which left-invariant
vector fields are parallel. Let ${\cal E}_I$ be the finite-dimensional vector
space of differential forms on $N_I$ which are parallel with
respect to $\nabla^{aff}$. 
Let $P : \Omega^*(N_I) \rightarrow \Omega^*(N_I)$ be orthogonal
projection onto ${\cal E}_I$, using $h(s)$.
From \cite[Proposition 1]{Lott (1999)}, $P$ is actually independent
of $s$ and arises from an averaging procedure over the group ${\frak N}_I$.
Let $\widehat{d}$ denote the exterior derivative on $\Omega^*(N_I)$, let
$\widehat{d}^*$ denote its adjoint with respect to $h(s)$ and put
\begin{equation} \label{Laplace}
\widehat{\triangle} \: = \: \widehat{d} \: \widehat{d}^* \: + \:
\widehat{d}^* \: \widehat{d}.
\end{equation}
The operators $\widehat{d}$, $\widehat{d}^*$ and $\widehat{\triangle}$ 
are diagonal with respect to the decomposition
\begin{equation} \label{decomp}
\Omega^*(N_I) \: = \: {\cal E}_I \: \oplus \: {\cal E}_I^\perp.
\end{equation}
We extend $\widehat{d}$, 
$\widehat{d}^*$ and $\widehat{\triangle}$ to act on $\Omega^*(N_I) \oplus
(ds \wedge \Omega^*(N_I))$, separately in each factor. 

Let $\{x^i\}$ be local coordinates on $N_I$.
Let $E^i$ denote exterior multiplication by $dx^i$ and let
$I_i$ denote interior multiplication by $\partial_{x^i}$.

\begin{lemma}
One has
\begin{equation} \label{V}
\partial_s \widehat{d}^* \: = \: \left[ \widehat{d}^*, \: V \right],
\end{equation}
where
\begin{equation} \label{VV}
V \: = \: 2 \: \sum_{i,j} S_i^{\: j} \: E^i \: I_j \: - \: \sum_i S^i_{\: i}.
\end{equation}
\end{lemma}
\begin{pf}
With our conventions, 
$\partial_s h \: = \: - \: 2 \: S$.
Given $\omega, \eta \in \Omega^*(N_I)$, let $\langle \omega, \eta \rangle
 \: \in \: 
C^\infty(N_I)$ be the inner product constructed using $h(s)$. One can check
that 
$\partial_s \: \langle \omega, \eta \rangle \: = \: 
\langle X \: \omega, \:  \eta \rangle$, where
$X \: = \: 2 \: \sum_{i,j} S_i^{\: j} \: E^i \: I_j$.
In addition, the derivative of the volume form is given by
$\partial_s \: d\vol \: = \: Y \: d\vol$, where
$Y \: = \: - \: \sum_i S^i_{\: i}$.
Differentiating the equation
\begin{equation}
\int_{N_I} \langle \widehat{d}^* \: \omega, \eta \rangle \: d\vol \: = \:
\int_{N_I} \langle \omega, \widehat{d} \: \eta \rangle \: d\vol 
\end{equation}
with respect to $s$ gives
\begin{align}
\int_{N_I} \langle X \: \widehat{d}^* \: \omega, \eta \rangle \: d\vol \: + \:
& \int_{N_I} \langle \partial_s \widehat{d}^* \: \omega, \eta \rangle 
\: d\vol \: + \:
\int_{N_I} \langle Y \: \widehat{d}^* \: \omega, \eta \rangle \: d\vol \: = \\
& \int_{N_I} \langle X \: \omega, \widehat{d} \: \eta \rangle \: d\vol \: + \:
\int_{N_I} \langle Y \: \omega, \widehat{d} \: \eta \rangle \: d\vol. \notag
\end{align}
As $\omega$ and $\eta$ are arbitrary, it follows that
\begin{equation}
X \: \widehat{d}^* \: + \: \partial_s \widehat{d}^* \: + \: Y \: 
\widehat{d}^* \: = \: \widehat{d}^* \: X \: + \: \widehat{d}^* \: Y,
\end{equation}
or
\begin{equation}
\partial_s \widehat{d}^* \: = \: 
[\widehat{d}^*,  X \: + \: Y].
\end{equation}
The lemma follows.
\end{pf}

Here $V$ is also diagonal 
with respect to the decomposition (\ref{decomp}).

It follows from Malcev's theorem that the harmonic forms on $(N_I, h(s))$
are parallel with respect to $\nabla^{aff}$.  
In particular, $\widehat{\triangle}$ is invertible
on ${\cal E}_I^\perp$. (Here ${\cal E}_I^\perp$ is also independent of $s$.) 
Let $G$ denote the corresponding Green's operator, which is the
inverse of $\widehat{\triangle}$
on ${\cal E}_I^\perp$ and which vanishes on ${\cal E}_I$.

\begin{lemma}
One has
\begin{equation}
\partial_s (\widehat{d}^* \: G) \: = \: - \: 
\left[ \widehat{d}, \: G \: \widehat{d}^* \: V \: \widehat{d}^* \: G \right].
\end{equation}   
\end{lemma}
\begin{pf}
Differentiating
the equations 
\begin{equation}
\widehat{\triangle} \: G \: = \: G \: \widehat{\triangle} \: = \: 1 \: - \: P
\end{equation}
and
\begin{equation}
P \: G \: = \: G \: P \: = \: 0
\end{equation}
with respect to $s$ gives
\begin{equation}
\partial_s G \: = \: - \: G \: (\partial_s \widehat{\triangle}) \: G.
\end{equation}
From (\ref{Laplace}),
\begin{equation}
\partial_s \widehat{\triangle} \: = \: \widehat{d} \: 
(\partial_s \widehat{d}^*) \: + \:
(\partial_s \widehat{d}^*) \: \widehat{d}.
\end{equation}
Then
\begin{align}
\partial_s (\widehat{d}^* \: G) \:  = \: & 
[\widehat{d}^*, V] \: G \: - \: \widehat{d}^* \: G \:
\left( \widehat{d} \: [\widehat{d}^*, V] \: + \: [\widehat{d}^*, V] \: 
\widehat{d} \right) G \: \\
= \: &
\widehat{d}^* \: V \: G \: - \:
V \: \widehat{d}^* \: G \: \notag \\
& - 
\: \widehat{d}^* \: G \: \widehat{d} \: \widehat{d}^* \: V \: G \: +
\: \widehat{d}^* \: G \: \widehat{d} \: V \widehat{d}^* \: G \: -
\: \widehat{d}^* \: G \: \widehat{d}^* \: V \: \widehat{d} \: G \: +
\: \widehat{d}^* \: G \: V \: \widehat{d}^* \: \widehat{d} \: G \notag \\ 
= \: &  \widehat{d}^* \: 
\left( I \: - \: 
\: G \: \widehat{d} \: \widehat{d}^* \right) \: V \: G 
- \: \left( I \: - \: 
\: \widehat{d}^* \: G \: \widehat{d} \right) \: V 
\: \widehat{d}^* \: G  \: +
\: \widehat{d}^* \: G \: V \: \widehat{d}^* \: \widehat{d} \: G \notag \\
\: = \: & \widehat{d}^* \: 
\left( I \: - \: 
\: G \: \widehat{\triangle} \right) \: V \: G 
- \: \left( I \: - \: G \: \widehat{\triangle} \: + \:
\: \widehat{d} \: G \: \widehat{d}^* \right) \: V 
\: \widehat{d}^* \: G  \: +
\: \widehat{d}^* \: G \: V \: \widehat{d}^* \: \widehat{d} \: G \notag \\
\: = \: & \widehat{d}^* \: 
P \: V \: G \: - \:  P \: V \: \widehat{d}^* \: G 
- \: \widehat{d} \: G \: \widehat{d}^*  \: V \: \widehat{d}^* \: G  \: +
\: \widehat{d}^* \: G \: V \: \widehat{d}^* \: \widehat{d} \: G \notag \\
\: = \: & \widehat{d}^* \: 
V \: P \: G \: - \:  V \: P \: \widehat{d}^* \: G 
- \: \widehat{d} \: G \: \widehat{d}^*  \: V \: \widehat{d}^* \: G  \: +
\: G \: \widehat{d}^* \: V \: \widehat{d}^* \: G \: \widehat{d} \notag \\
\: = \: &
- \: [\widehat{d}, \: G \: \widehat{d}^*  \: V \: 
\widehat{d}^* \: G]. \notag
\end{align}
This proves the lemma.
\end{pf}

Let $e(ds)$ denote exterior multiplication by $ds$.
Define ${\cal K} \: : \: \Omega^*(\overline{U_I}) \rightarrow 
\Omega^*(\overline{U_I})$ by
\begin{equation}
{\cal K} \: = \: \widehat{d}^* \: G \: - \: e(ds) \:  G \: \widehat{d}^* \:
V \: \widehat{d}^* \: G. 
\end{equation}

\begin{lemma}
Acting on $\Omega^*(\overline{U_I})$, one has
\begin{equation} \label{homotopy}
d \: {\cal K} \: + \: {\cal K} \: d \: = \: 1 \: - \: P.
\end{equation}
(In this last equation, $P$ acts fiberwise.)
\end{lemma}
\begin{pf}
Using the fact that
\begin{equation}
d \: = \: \widehat{d} \: + \: e(ds) \: \partial_s,
\end{equation}
we have
\begin{align}
d \: {\cal K} \: + \: {\cal K} \: d \: = \: &
\left(\widehat{d} \: + \: e(ds) \: \partial_s \right) 
\left(\widehat{d}^* \: G \: - \: e(ds) \:  G \: \widehat{d}^* \:
V \: \widehat{d}^* \: G \right) \: + \\
& \left( \widehat{d}^* \: G \: - \: e(ds) \:  G \: \widehat{d}^* \:
V \: \widehat{d}^* \: G \right) 
\left( \widehat{d} \: + \: e(ds) \: \partial_s\right) \notag \\
= \: & \: \widehat{d} \: \widehat{d}^* \: G \: + \: \widehat{d}^* \: G \:
\widehat{d} \: + \: e(ds) \: \left( [\partial_s, \: \widehat{d}^* \: G] \:
+ \: [\widehat{d}, \: G \: \widehat{d}^* \: V \: \widehat{d}^* \: G] \right)
\notag \\
= \: & \: I \: - \: P. \notag
\end{align}
This proves the lemma.
\end{pf}

Consider the trivial vector bundle 
${\cal W}_I \: = \: [0, \infty) \times {\cal E}_I$ over
$[0, \infty)$. Let $\widehat{d}^{inv}$ be the restriction of
$\widehat{d}$ to ${\cal E}_I \subset \Omega^*(N_I)$ and consider the
flat superconnection $A_I$ on ${\cal W}_I$ whose action on 
$\Omega^*([0, \infty); {\cal W}_I)$ is given by
\begin{equation}
A_I \: = \: \widehat{d}^{inv} \: + \: e(ds) \: \partial_s.
\end{equation}
That is, $A_I$ is simply the restriction of $d$ from
$\Omega^*(\overline{U_I})$ to $\Omega^*([0, \infty); {\cal W}_I)$. Then
(\ref{homotopy}) gives a homotopy equivalence between the cochain complexes
$\left( \Omega^*(\overline{U_I}), d \right)$ and
$\left( \Omega^*([0, \infty); {\cal W}_I), A_I \right)$.  

From the Gauss-Codazzi equation and the results of 
\cite{Cheeger-Fukaya-Gromov (1992)}, we can assume that there is 
a uniform upper bound on the absolute values of the sectional
curvatures of the fibers $(N_I, h(s))$, of the form $\const \: b^2$. Then from
\cite[Proposition 2]{Lott (1999)}, it follows that
if $\delta$ is small enough then there
is a uniform positive lower bound on the eigenvalues of 
$\widehat{\triangle} \big|_{{\cal E}^\perp}$. Hence
${\cal K}$ is a bounded operator.  
Then it follows as in \cite[Lemma 1]{Lott (1997)}
that the (reduced) $L^2$-cohomology of $\overline{U_I}$ is isomorphic
to the (reduced) $L^2$-cohomology of 
$\left( \Omega^*([0, \infty); {\cal W}_I), A_I \right)$, 
where $\Omega^*([0, \infty); {\cal W}_I)$ 
acquires an $L^2$-inner product from $\Omega^*(\overline{U_I})$.
From Hodge theory, the (reduced) $L^2$-cohomology of 
$\left( \Omega^*([0, \infty); {\cal W}_I), A_I \right)$
is isomorphic to the vector space of square-integrable solutions
to the equation
\begin{equation} \label{solns}
(A_I A_I^* \: + \:  A_I^* A_I) \: \psi \: = \: 0
\end{equation}
on $[0, \infty)$, with absolute boundary
conditions at $\{0\}$. 
However, as $A_I A_I^* \: + \:  A_I^* A_I$ is a second-order 
ordinary differential operator, the
solution space of (\ref{solns}) is finite-dimensional.  
This proves the theorem.

\section{Geometry of Finite-Volume Negatively-Curved manifolds} \label{sect3}

We review some results from \cite{Eberlein (1980)} and 
\cite{Heintze-Im Hof (1977)}.
Let $(M,g)$ 
be a complete connected Riemannian manifold of finite volume
whose sectional curvatures satisfy $-b^2 \le K \le - a^2$, with $0 < a \le b$. 
Then $M$ is diffeomorphic to the interior of a smooth compact connected
manifold-with-boundary $\overline{M}$. 
The boundary components of $\overline{M}$ are diffeomorphic to
infranilmanifolds. If $N$ is such a
boundary component then there is a corresponding end $E$
of $M$. Let $s$ be a Busemann function for a ray exiting $E$. 
Then after changing $s$ by a constant if necessary, 
there are a neighborhood $U$ of $E$ and a 
$C^1$-diffeomorphism ${\cal F} : (0, \infty) \times N \rightarrow U$
so that 
\begin{equation}
{\cal F}^* \left( g \big|_U \right) = ds^2 + h(s),
\end{equation}
where for $s \in (0, \infty)$,  $h(s)$ is a Riemannian metric on $N$.
We will think of $s$ as a coordinate function on $U$.
The slices $N(s) = \{s\} \times N$ are projections of horospheres in the
universal cover $\widetilde{M}$.  {\it A priori}, the Busemann function is
only $C^2$-smooth on $M$ and the Riemannian metric $h(s)$ is only 
$C^1$-smooth on $N$. 
Given $n \in N$, the curve $s \rightarrow (s, n)$ is a unit-speed geodesic 
which intersects the slices orthogonally.
All of the rays in $M$ which exit $E$ arise in this way. 

As $s$ is $C^2$-smooth, 
the second fundamental form $S(s)$ of $N(s)$ exists and is
continuous on $N(s)$.
From Jacobi field estimates, it satisfies
\begin{equation}
a \: h(s) \: \le \: S(s) \: \le \: b \: h(s)
\end{equation}
and the metric $h(s)$ satisfies
\begin{equation}
e^{-2bs} h(0) \: \le \: h(s) \: \le \: e^{-2as} h(0).
\end{equation}

\section{Infranilmanifolds} \label{sect4}

Let $N$ be a boundary component of $\overline{M}$. It has a regular
covering by a nilmanifold $\Gamma \backslash {\frak N}$, with covering group
$F$. Here ${\frak N}$ is a simply-connected connected nilpotent Lie group,
$\Gamma$ is a lattice in ${\frak N}$ and $F$ is a finite group of
automorphisms of ${\frak N}$ which preserve $\Gamma$.
Let ${\frak n}$ be the Lie algebra of ${\frak N}$. 
Let $\Lambda^*({\frak n}^*)^F$ denote
the $F$-invariant subspace of $\Lambda^*({\frak n}^*)$. Let $\triangle_{N(s)}$
denote the differential form Laplacian on $N(s)$ (which can be defined using
quadratic forms \cite[Vol. I, Theorem VIII.15]{Reed-Simon (1978)} 
even if $h(s)$ is only $C^1$-smooth). Given 
$\lambda \in [0, \infty)$, let $P_{N(s)}(\lambda)$ denote the spectral
projection onto the direct sum of the eigenspaces of $\triangle_{N(s)}$ with
eigenvalue less than or equal to $\lambda$.

\begin{proposition} \label{prop1}
There are constants $c_1, c_2 > 0$ 
such that for all sufficiently large $s$, the
images of 
$P_{N(s)}(c_1^2 \: b^2)$ and $P_{N(s)}(c_2^2 \: a^2 \:  e^{2as})$ are 
isomorphic to $\Lambda^*({\frak n}^*)^F$.
\end{proposition}
\begin{pf}
Suppose first that the parametrization
${\cal F} : (0, \infty) \times N \rightarrow U$
is smooth. From the Gauss-Codazzi equation, the
intrinsic sectional curvatures $R^{N(s)}$
of $N(s)$ are bounded in absolute value
by a universal constant times $b^2$. From 
\cite{Cheeger-Fukaya-Gromov (1992)}, there is an $\epsilon > 0$ such that
for all $s \in [1, \infty)$, there is a metric $h_0(s)$ on $N(s)$, 
coming from an $F$-invariant left-invariant inner product on ${\frak N}$, with
\begin{equation}
e^{-\epsilon} h_0(s) \: \le \: h(s) \: \le \: e^{\epsilon} h_0(s).
\end{equation} 
By \cite{Dodziuk (1982)}, there is an integer $J > 0$ such that
the $j$-th eigenvalue $\lambda_{p,j}$ of the 
$p$-form Laplacian satisfies
\begin{equation} \label{cont}
e^{-J \epsilon} \lambda_{p,j}(h_0(s)) \: \le \: \lambda_{p,j}(h(s)) \: \le \: 
e^{J \epsilon} \lambda_{p,j}(h_0(s)).
\end{equation} 
Thus without loss of generality, we may assume that $h(s)$ comes from 
a left-invariant inner product on ${\frak N}$.

The vector space of
$F$-invariant left-invariant differential forms on ${\frak N}$ 
is isomorphic to $\Lambda^*({\frak n}^*)^F$. These differential forms push
down to comprise a vector space ${\cal V}$ of
differential forms on $N(s)$. The Laplacian $\triangle_{N(s)}$ has an
orthogonal direct sum decomposition
\begin{equation}
\triangle_{N(s)} = \triangle_{\cal V} \oplus 
\triangle_{{\cal V}^\perp}.
\end{equation}
From \cite[Proposition 2]{Lott (1999)}, there is a constant $c_2 > 0$ such
for sufficiently large $s$, 
the eigenvalues of $\triangle_{{\cal V}^\perp}$ are greater than 
$c_2^2 \: a^2 \:  e^{2as}$.

It remains to show that there is a constant $c_1 > 0$ such that
the eigenvalues of $\triangle_{\cal V}$ 
are less than or equal to $c_1^2 \: b^2$, uniformly
in $s$. We follow the notation of \cite[Section 3]{Lott (1999)}.
Let $\{e_i\}$ be the orthonormal basis of ${\frak n}$ described in 
\cite[Section 3]{Lott (1999)}, with dual
basis $\{\tau^i\}$. Let $E^i$ denote exterior multiplication by $\tau^i$ and
let $I_i$ denote interior multiplication by $e_i$.
The exterior derivative $d$, acting on $\Omega^*(N(s))$, can be written as
$d = \sum_i E^i \: \nabla^{N(s)}_{e_i}$, and its adjoint can be written as
$d^* = - \sum_i I^i \: \nabla^{N(s)}_{e_i}$. Now $\nabla^{N(s)}_{e_i}$ acts on 
${\cal V}$ as $\sum_{j,k} \omega^j_{\: ki} \: E^j \: I^k$, where
$\{ \omega^j_{\: ki} \}$ are the components of the Levi-Civita connection
$1$-form $\omega$ of the left-invariant metric. 
Hence $\triangle_{\cal V}$ is quadratic in $\omega$.
From \cite[Lemma 3]{Lott (1999)}, there is a constant,
which only depends on $\dim(N)$, such that
\begin{equation}
\parallel \omega \parallel_\infty^2 \: \le \: \const \: \parallel
R^{N(s)} \parallel_\infty.
\end{equation}
The proposition follows, under the assumption that the parametrization 
${\cal F} : (0, \infty) \times N \rightarrow U$ is smooth.

In the general case,
thinking of $N(s)$ as the graph of a $C^2$-function on $N$, for any 
$\epsilon > 0$ we can
find a smooth hypersurface $N^\prime$ of $M$ which is $\epsilon$-close to
$N(s)$ in the $C^2$-topology. Then the proposition holds for $N^\prime$.
Using the continuity of the eigenvalues with respect to the metric,
in the $C^0$-topology, as in
(\ref{cont}), the proposition follows. In fact, we can take $c_1$ and $c_2$
to only depend on $\dim(N)$, although we will not need this.
\end{pf}

Let $\{N_I\}_{I=1}^B$ be the boundary components of
$\overline{M}$, with corresponding ends $E_I$ and neighborhoods of the ends
$U_I$. By reducing $U_I$ if necessary, we may assume that
Proposition \ref{prop1} holds for all $s \: \ge \: 0$, with $c_1 b \: < \: 
c_2 a$.
As in \cite[Proposition 2.1]{Donnelly-Li (1979)}, 
the essential spectrum of $\triangle^M_p$ is invariant under
compactly-supported changes of the metric.
Thus without loss of generality,
we may assume that the metric on $\overline{U_I}$ is a product near
$\{0\} \times N_I$, with
Proposition \ref{prop1} still holding for $s \: \ge \: 0$.
Let $\Omega^*_{I}$
denote the smooth compactly-supported forms 
on $[0, \infty) \times N_I$ which satisfy relative boundary conditions
at $\{0\} \times N_I$.
Let $H^\prime$ be the $L^2$-completion of $\bigoplus_{I=1}^B 
\Omega^*_{I}$.
The Laplacian $\triangle^\prime \: = \: d \: d^* \: + \: d^* \: d$, 
defined initially on $\bigoplus_{I=1}^B \Omega^*_{I}$, is
a densely-defined self-adjoint operator on $H^\prime$ and corresponds
to relative boundary conditions.

For later use, we write $d$ and $d^*$ more explicitly. Fix $I$.
Let $\{x^i\}_{i=1}^{n-1}$ be local coordinates on $N_I$ and write the metric on
$U_I$ as
$ds^2 \: + \: \sum_{i,j} h_{ij} \:  dx^i \: dx^j$. We think
of $s = x^0$ as another coordinate.
Let $S_{ij}$ be the second fundamental form of $\{s\} \times N_I$.
We let Greek letters run over $\{0, \ldots, n - 1\}$ and we let Roman letters
run over $\{1, \ldots, n - 1\}$.
The nonzero Christoffel symbols are
\begin{align}
\Gamma_{ijk} & = \Gamma_{ijk}(h), \\
\Gamma_{0ij} & = S_{ij}, \notag \\
\Gamma_{ij0} & = - \: S_{ij}, \notag \\
\Gamma_{i0j} & = - \: S_{ij}. \notag
\end{align}
Let $E^\alpha$ denote exterior multiplication by $dx^\alpha$ and let
$I_\alpha$ denote interior multiplication by $\partial_{x^\alpha}$.
Covariant differentiation on forms is given in local coordinates by
\begin{equation}
\nabla_{\partial_{x^\alpha}} \: = \: \partial_{x^\alpha} \: - \:
\sum_{\beta, \gamma} \Gamma^{\gamma}_{\: \beta \alpha} \: E^\beta \: 
I_\gamma.
\end{equation}
Let $\widehat{\nabla}$ denote the covariant derivative on $N_I(s)$.
Then
\begin{align} \label{covder}
\nabla_{\partial_{x^i}} \: & = \: \widehat{\nabla}_{\partial_{x^i}} \: 
 - \:
\sum_{j} S_{ij} \: E^j \: 
I_0 + \:
\sum_{j} S_i^{\: j} \: E^0 \: I_j, \\
\nabla_{\partial_{s}} \: & = \: \partial_{s} \: + \:
\sum_{i,j} S^i_{\: j} \: E^j \: I_i. \notag
\end{align}
Let $\widehat{d} \: = \: \sum_i E^i \: \widehat{\nabla}_{\partial_{x^i}}$ be
the exterior derivative on $N_I(s)$, extended to act on
$\Omega^*(N_I(s)) \oplus (ds \wedge \Omega^*(N_I(s)))$, and let 
$\widehat{d}^* \: = \: - \: \sum_i I^i \: \widehat{\nabla}_{\partial_{x^i}}$ be
its adjoint.
Then
\begin{equation} \label{explicit}
d \:  = \: \widehat{d} \: + 
E^0 \: \partial_{s}
\end{equation}
and
\begin{align} \label{explicit2}
d^*  \:  & = \:  - \: \sum_\alpha I^\alpha \: 
\nabla_{\partial_{x^\alpha}} \\
& = \: \widehat{d}^* \: - \:
I^0 \: \left( \partial_{s} \: + \: \sum_{i,j} 
S^i_{\: j} \: (E^j \: I_i \: - I_i \: E^j) \right). \notag
\end{align}

\section{Boundedness of the Off-Diagonal Operators} \label{sect6}

Given $I$, consider $N_I$ to be an infranilmanifold which is $F_I$-covered by a
nilmanifold $\Gamma_I \backslash {\frak N}_I$ and let ${\frak n}_I$ be the
Lie algebra of ${\frak N}_I$.
Let $E_I = [0, \infty) \times \Lambda^*({\frak n}_I^*)^{F_I}$ be the trivial
vector bundle on $[0, \infty)$ with fiber $\Lambda^*({\frak n}_I^*)^{F_I}$.

Let $\Omega^*([0, \infty); E_I)$ be the smooth compactly-supported
forms on $[0, \infty)$, with value in $E_I$.
Using Proposition \ref{prop1}, there 
is an embedding of $\Omega^*([0, \infty); E_I)$ into
$\Omega^*(\overline{U_I})$. Let $\Omega^*_{rel}([0, \infty); E_I)$ be the
subspace of $\Omega^*([0, \infty); E_I)$ consisting of elements which
satisfy relative
boundary conditions at $\{0\}$.
Let $H_0$ be the completion of 
$\bigoplus_{I=1}^B \Omega^*_{rel}([0, \infty); E_I)$ in $H^\prime$ and let 
$H_1$ be its orthogonal complement. Roughly speaking, the elements of
$H_0$ correspond to fiberwise low-energy forms
and the elements of $H_1$ correspond to
fiberwise high-energy forms.

Let $P_0 : H^\prime \rightarrow H_0$ and $P_1 : H^\prime \rightarrow H_1$ be
the orthogonal projections.
With respect to the orthogonal decomposition $H^\prime \: = \: H_0 \oplus
H_1$, write
\begin{equation}
d \: = \:
\begin{pmatrix}
{\cal A} & {\cal B} \\
{\cal C} & {\cal D}
\end{pmatrix}.
\end{equation}
Then
\begin{equation}
d^* \: = \:
\begin{pmatrix}
{\cal A}^* & {\cal C}^* \\
{\cal B}^* & {\cal D}^*
\end{pmatrix}
\end{equation}
and
\begin{equation}
\triangle^\prime \: = \:
\begin{pmatrix}
{\cal A} {\cal A}^* \: + \: {\cal A}^* {\cal A} \: + \:
{\cal B} {\cal B}^* \: + \: {\cal C}^* {\cal C} & 
{\cal A} {\cal C}^* \: + \: {\cal B} {\cal D}^* \: + \: 
{\cal A}^* {\cal B} \: + \: {\cal C}^* {\cal D} \\
{\cal C} {\cal A}^* \: + \:{\cal D} {\cal B}^* \: + \: 
{\cal B}^* {\cal A} \: + \: {\cal D}^* {\cal C} & 
{\cal D} {\cal D}^* \: + \: {\cal D}^* {\cal D} \: + \: 
{\cal B}^* {\cal B} \: + \: {\cal C} {\cal C}^*
\end{pmatrix}.
\end{equation}

\begin{proposition} \label{prop3}
The operators ${\cal B} \: : \: H_1
\rightarrow H_0$ and ${\cal C} \: : \: H_0
\rightarrow H_1$ are bounded.
\end{proposition}
\begin{pf}
We have
\begin{equation}
{\cal B} \: = \: P_0 \: d \: P_1 \: = \: \left( P_1 \: d^* \: P_0 \right)^*
\end{equation}
and
\begin{equation}
{\cal C} \: = \: P_1 \: d \: P_0
\end{equation}
From (\ref{explicit}) and (\ref{explicit2}), in order 
to show that ${\cal B}$ and ${\cal C}$ are bounded 
it is enough to show that 
\begin{equation}
P_1 \: \partial_s \: P_0 \: = \: P_1 \: (\partial_s P_0)
\end{equation}
is bounded.

Let $\gamma$ be the circle of radius $c_1 \: b$ 
around the origin in $\C$, oriented
counterclockwise.  From Proposition \ref{prop1},
\begin{equation} \label{projection}
P_0(s) \: = \: \oint_\gamma (\lambda \: - \: \widehat{d} \: - 
\widehat{d}^*)^{-1} \: \frac{d\lambda}{2 \pi i}.
\end{equation}
Here $P_0(s)$ is a projection on $\bigoplus_{I=1}^B \left( \Omega^*(N_I(s))
\: \oplus \: (ds \:\wedge \: \Omega^*(N_I(s))) \right)$.
We note that the Hilbert space structure on 
$\bigoplus_{I=1}^B \left( \Omega^*(N_I(s))
\: \oplus \: (ds \:\wedge \: \Omega^*(N_I(s))) \right)$ depends on $s$, but the
underlying topological vector space structure on
$\bigoplus_{I=1}^B \left( \Omega^*(N_I)
\: \oplus \: (ds \:\wedge \: \Omega^*(N_I)) \right)$
does not.  Hence it makes sense
to differentiate (\ref{projection}) with respect to $s$, giving
\begin{align}  \label{diff}
\partial_s P_0 \: = \: & \oint_\gamma (\lambda \: - \: \widehat{d} \: - 
\widehat{d}^*)^{-1} \: \partial_s (\widehat{d} \: + \: \widehat{d}^*) \: 
(\lambda \: - \: \widehat{d} \: - 
\widehat{d}^*)^{-1} \: \frac{d\lambda}{2 \pi i} \\
= \: & \oint_\gamma (\lambda \: - \: \widehat{d} \: - 
\widehat{d}^*)^{-1} \: \partial_s \widehat{d}^* \: 
(\lambda \: - \: \widehat{d} \: - 
\widehat{d}^*)^{-1} \: \frac{d\lambda}{2 \pi i} \notag \\
= \: & \oint_\gamma (\lambda \: - \: \widehat{d} \: - 
\widehat{d}^*)^{-1} \: [\widehat{d}^*, V] \:
(\lambda \: - \: \widehat{d} \: - 
\widehat{d}^*)^{-1} \: \frac{d\lambda}{2 \pi i}, \notag
\end{align}
where $V$ is as in (\ref{VV}).

As $P_1(s) \: = \: 1 \: - \: P_0(s)$, it 
follows from differentiating $P_0^2(s) \: = \: P_0(s)$ that
\begin{equation}
P_1 \: (\partial_s P_0) \: = \: (\partial_s P_0) \: P_0.
\end{equation}
If $\eta_0 \: \in \: \Image(P_0)$ is an eigenform for $\widehat{d} \: + \:
\widehat{d}^*$ with eigenvalue $\lambda_0$ and
$\eta_1 \: \in \: \Image(P_1)$ is an eigenform for $\widehat{d} \: + \:
\widehat{d}^*$ with eigenvalue $\lambda_1$ then
\begin{align}
\langle \eta_1, \: \oint_\gamma (\lambda \: - \: \widehat{d} \: - 
\widehat{d}^*)^{-1} \: & [\widehat{d}^*, V] \:
(\lambda \: - \: \widehat{d} \: - 
\widehat{d}^*)^{-1} \: \frac{d\lambda}{2 \pi i} \: \eta_0 \rangle \\
& = \: \oint_\gamma \langle \eta_1, \: (\lambda \: - \: \lambda_1)^{-1} \: 
[\widehat{d}^*, V] \:
(\lambda \: - \: \lambda_0)^{-1} 
\: \eta_0 \rangle \: \frac{d\lambda}{2 \pi i} \notag \\
& = \: - \: 
\frac{1}{\lambda_1 \: - \: \lambda_0} \:
\langle \eta_1, \: [\widehat{d}^*, V] \: \eta_0 \rangle. \notag
\end{align}
It follows that
\begin{align} \label{follows}
\oint_\gamma (\lambda \: - \: \widehat{d} \: - 
\widehat{d}^*)^{-1} \: & [\widehat{d}^*, V] \:
(\lambda \: - \: \widehat{d} \: - 
\widehat{d}^*)^{-1} \: \frac{d\lambda}{2 \pi i} \: \eta_0 \\
& = \: - \: 
\left( (\widehat{d} \: + \widehat{d}^*) \Big|_{Im(P_1(s))} \: - \:
\lambda_0 \right)^{-1} 
P_1(s) \:  [\widehat{d}^*, V] \:
\eta_0. \notag
\end{align}

Using (\ref{diff}) and (\ref{follows}),
in order to prove the proposition it suffices to show that as
$\eta_0$ runs over unit-length eigenforms of 
$(\widehat{d} \: + \: \widehat{d}^*) \Big|_{Im(P_0(s))}$,
one has a bound on the norm of 
$$\left(
(\widehat{d} \: + \widehat{d}^*) \Big|_{Im(P_1(s)))} \: - \lambda_0
\right)^{-1} 
P_1(s) \:  [\widehat{d}^*, V] \:
\eta_0$$ which is uniform in $s$. Writing
\begin{align}
\left((\widehat{d} \: + \widehat{d}^*) \Big|_{Im(P_1(s)))} \: - \lambda_0
\right)^{-1} 
& P_1(s) \:  [\widehat{d}^*, V] \:
\eta_0 \: = \\
& \left((\widehat{d} \: + \widehat{d}^*) \Big|_{Im(P_1(s)))} \: - 
\lambda_0 \right)^{-1} 
P_1(s) \:  \widehat{d}^* \: (V \: \eta_0) \: - \notag \\
& \left( (\widehat{d} \: + \widehat{d}^*) \Big|_{Im(P_1(s)))} \: - 
\lambda_0 \right)^{-1} 
P_1(s) \: V \:  \widehat{d}^* \eta_0, \notag
\end{align}
we know from (\ref{VV}) and Proposition \ref{prop1}
that we have bounds on $\Big| V \: \eta_0 \Big|$ and
$\Big| V \: \widehat{d}^* \eta_0 \Big|$ given by
$\const \: b$ and $\const \: b^2$, respectively. Hence it suffices to
show that the operators 
\begin{equation} \label{op1}
\left(
(\widehat{d} \: + \widehat{d}^*) \Big|_{Im(P_1(s)))} \: - \lambda_0
\right)^{-1} 
\:  \widehat{d}^* \: P_1(s)
\end{equation}
and 
\begin{equation} \label{op2}
\left( (\widehat{d} \: + \widehat{d}^*) \Big|_{Im(P_1(s)))} \: - 
\lambda_0 \right)^{-1} 
P_1(s)
\end{equation}
have uniform bounds on their operator norms.

Put $\widehat{\triangle} \: = \: (\widehat{d} \: + \: \widehat{d}^*)^2$.
As (\ref{op1}) vanishes on $\Image ( \widehat{d}^* )$, using the
Hodge decomposition it is enough to consider its action on 
$\Image(\widehat{d})$. Given $\eta \in \Image(P_1(s))$, we have
\begin{align}
\left( (\widehat{d} \: + \widehat{d}^*) \: - \: \lambda_0 \right)^{-1} 
\:  \widehat{d}^* \: \widehat{d} \eta \: = \: &
\frac{\widehat{d} \: + \widehat{d}^* \: + \: \lambda_0}{\widehat{\triangle}
\: - \: \lambda_0^2} \:  \widehat{d}^* \: \widehat{d} \eta \\
= \: & \frac{\widehat{\triangle}}{\widehat{\triangle} - \lambda_0^2} \:
\widehat{d} \eta \: 
+ \: \frac{\lambda_0}{\widehat{\triangle} - \lambda_0^2} \:  
\widehat{d}^* \: \widehat{d} \eta. \notag
\end{align}
By Proposition \ref{prop1}, the operator norm of 
$\frac{\widehat{\triangle}}{\widehat{\triangle} - \lambda_0^2}$, acting
on $\Image(P_1(s))$, is at most 
$\frac{c_2^2 \: a^2 \: e^{2as}}{c_2^2 \: a^2 \: e^{2as} \: - \: 
c_1^2 \: b^2}$. If 
${\cal T} \: = \: \frac{\lambda_0}{\widehat{\triangle} - \lambda_0^2} \:  
\widehat{d}^*$ then 
${\cal T}^* \: {\cal T} \: = \: \frac{\lambda_0^2 \: \widehat{d} \:
\widehat{d}^*}{(\widehat{\triangle} - \lambda_0^2)^2}
$
which, acting on $\Image(\widehat{d})$, equals
$
\frac{\lambda_0^2 \: 
\widehat{\triangle}}{(\widehat{\triangle} - \lambda_0^2)^2}
$.
Hence the norm of ${\cal T}$, 
acting on $\Image(\widehat{d}) \bigcap \Image(P_1(s))$,
is at most
$\frac{c_1 \: b \: c_2 \: a \: e^{as}}{c_2^2 \: a^2 \: e^{2as} \: - \:
c_1^2 \: b^2}$.

By Proposition \ref{prop1},
the operator norm of (\ref{op2}) is at most
$\frac{1}{c_2 \: a \: e^{as} \: - \: c_1 \: b}$. 
The proposition follows.
\end{pf}

\section{High Energy Forms} \label{sect5}

\begin{proposition} \label{prop2}
The operator ${\cal D} {\cal D}^* \: + \: 
{\cal D}^* {\cal D} \: + \: {\cal B}^* {\cal B} \: + \: 
{\cal C} {\cal C}^*$ 
has vanishing essential spectrum.
\end{proposition}
\begin{pf}
Without loss of generality, we consider the neighborhood $U_I$ of a single end.
By standard arguments \cite{Donnelly-Li (1979)}, 
it suffices to show that as $c \rightarrow \infty$,
the infimum of 
\begin{equation}
| {\cal D} J|^2 \: + \: | {\cal D}^* J|^2 \: + \: | {\cal B} J|^2 \: + \: 
| {\cal C}^* J|^2,
\end{equation}
where $J$ runs over smooth unit-length elements of
$H_1$ with compact support in $[c, \infty) \times N_I$, goes to infinity.
In this proof, all norms will be $L^2$-norms on 
$\overline{U_I} \: = \: [0, \infty) \: \times \: N_I$.

Taking $c \: > \:  0 $, we can ignore boundary terms in the following
equations. We have
\begin{equation}
| {\cal D} J|^2 \: + \: | {\cal D}^* J|^2 \:  = \:
|d J \: - \: {\cal B} J|^2 \: + \: |d^* J \: - \: {\cal C}^* J|^2.
\end{equation}
From Proposition \ref{prop3}, ${\cal B}$ and ${\cal C}^*$ are bounded. 
In terms of the two-component vector 
$(d J \: - \: {\cal B} J, d^* J \: - \: {\cal C}^* J)$, we can write
\begin{align}
|d J \: - \: {\cal B} J|^2 \: + \: |d^* J \: - \: {\cal C}^* J|^2 \: & =
\: |(d J \: - \: {\cal B} J, d^* J \: - \: {\cal C}^* J)|^2 \\
& = \: |(d J, d^* J) \: - \: 
({\cal B} J, {\cal C}^* J)|^2 \notag \\
& \ge \: 
\left( |(d J, d^* J)| \: - \: |({\cal B} J, {\cal C}^* J)| \right)^2
\notag \\
& \ge 
\left( \max( \sqrt{|d J|^2 \: + \: |d^* J|^2} \: - \: \const \: |J|, 0)
\right)^2,
\end{align}
where ``$\const$'' in this proof will denote something that is independent of 
$c$. 
Hence it suffices to consider $| d J|^2 \: + \: | d^* J|^2$.

From Bochner's equation,
\begin{equation}
| d J|^2 \: + \: | d^* J|^2 \: = \:
|\nabla J|^2 \: + \: \int_{U_I} \sum_{p,q,r,s = 0}^{n-1}
R^M_{pqrs} \langle E^p I^q J, E^r I^s J \rangle \: \ge \:
|\nabla J|^2 \: - \: \const \: |J|^2. 
\end{equation}
Letting $\nabla^{vert}$ denote covariant differentiation in vertical
directions, we have
\begin{equation}
|\nabla J|^2 \: = \: |\nabla^{vert} J|^2 \: + \: |\nabla_{\partial_s} J|^2.
\end{equation}
Thus 
\begin{equation}
| d J|^2 \: + \: | d^* J|^2 \: \ge \:
|\nabla^{vert} J|^2 \: - \: \const \: |J|^2 \: \ge \:
\left( |\widehat{\nabla} J| \: - \: 
|(\nabla^{vert}  - \widehat{\nabla}) J | \right)^2 \: - \: \const \: |J|^2. 
\end{equation}
Using (\ref{covder}), we obtain
\begin{equation}
| d J|^2 \: + \: | d^* J|^2 \: \ge \:
\max \left( |\widehat{\nabla} J| \: - \: \const \: |J|, 0 \right)^2 \: - \:
\const \: |J|^2. 
\end{equation}

Applying Bochner's equation fiberwise gives
\begin{equation}
| \widehat{d} J|^2 \: + \: | \widehat{d}^* J|^2 \: = \:
|\widehat{\nabla} J|^2 \: + \: \int_{U_I} \sum_{p,q,r,s = 1}^{n-1}
R^Z_{pqrs} \langle E^p I^q J, E^r I^s J \rangle.
\end{equation}
From the Gauss-Codazzi equation,
\begin{equation} \label{GC}
\int_{U_I} \sum_{p,q,r,s = 1}^{n-1}
R^Z_{pqrs} \langle E^p I^q J, E^r I^s J \rangle \: = \:
\int_{U_I} \sum_{p,q,r,s = 1}^{n-1}
(R^M_{pqrs} \: + \: S_{pr} \: S_{qs} \: - \: S_{ps} \: S_{qr}) \: 
\langle E^p I^q J, E^r I^s J \rangle.
\end{equation}
Hence 
\begin{equation} \label{GC2}
|\widehat{\nabla} J|^2 \: \ge \: 
| \widehat{d} J|^2 \: + \: | \widehat{d}^* J|^2 
 \: - \: \const \: |J|^2. 
\end{equation}
(We note that the right-hand-side of (\ref{GC}) makes sense
even if the Busemann function
$s$ is only $C^2$-smooth. Hence (\ref{GC2}) is valid in this generality.)
From Proposition \ref{prop1}, we have
\begin{equation}
|\widehat{d} J|^2 \: + \: |\widehat{d}^* J|^2 \: \ge \:
c_2^2 \: a^2 \: e^{2ac} \: |J|^2. 
\end{equation}
Taking $c \rightarrow \infty$, the proposition follows.
\end{pf}

\section{Proof of Theorem \ref{thm2}} \label{sect7}

We will use the general identity that
\begin{equation} \label{inverse}
\begin{pmatrix}
\alpha & \beta \\
\gamma & \delta
\end{pmatrix}^{-1} \: = \:
\begin{pmatrix}
(\alpha \: - \: \beta \: \delta^{-1} \: \gamma)^{-1} & 
- \: (\alpha \: - \: \beta \: \delta^{-1} \: \gamma)^{-1} \: \beta \:
\delta^{-1} \\
- \: \delta^{-1} \: \gamma \: 
(\alpha \: - \: \beta \: \delta^{-1} \: \gamma)^{-1} & 
\delta^{-1} \: + \: \delta^{-1} \: \gamma \:
(\alpha \: - \: \beta \: \delta^{-1} \: \gamma)^{-1} \: \beta \: \delta^{-1}
\end{pmatrix},
\end{equation}
provided that $\delta$ and 
$\alpha \: - \: \beta \: \delta^{-1} \: \gamma$ are invertible.

By a standard argument as in \cite[Proposition 2.1]{Donnelly-Li (1979)},
the essential spectra of $\triangle^M_p$ and
$\triangle^\prime_p$ are the same. For simplicity, we will omit the
subscript $p$. Using Proposition \ref{prop2}, 
it is enough to show that $\triangle^\prime$ and
\begin{equation}
{\cal L} \: = \: 
\begin{pmatrix}
{\cal A} {\cal A}^* \: + \: {\cal A}^* {\cal A} & 0 \\ 
0 & {\cal D} {\cal D}^* \: + \: {\cal D}^* {\cal D} \: + \: 
{\cal B}^* {\cal B} \: + \: {\cal C} {\cal C}^*
\end{pmatrix}.
\end{equation}
have the same essential spectra. To show this, from 
\cite[Vol. IV, Theorem XIII.14]{Reed-Simon (1978)} it
suffices to show that $\left( \triangle^\prime \: + \: k \: I \right)^{-1} \: -
\: \left( {\cal L} \: + \: k \: I \right)^{-1}$ is compact for some
$k \: > \: 0$.

We put
\begin{equation}
\begin{pmatrix}
\alpha & \beta \\
\gamma & \delta
\end{pmatrix} \: = \: \triangle^\prime \: + \: k \: I.
\end{equation}
Explicitly,
\begin{align}
\alpha \: & = \: 
{\cal A} {\cal A}^* \: + \: {\cal A}^* {\cal A} \: + \:
{\cal B} {\cal B}^* \: + \: {\cal C}^* {\cal C}  \: + \: k \: I, \\ 
\beta \: & = \: {\cal A} {\cal C}^* \: + \: {\cal B} {\cal D}^* \: + \: 
{\cal A}^* {\cal B} \: + \: {\cal C}^* {\cal D} \notag \\
\gamma \: & = \: {\cal C} {\cal A}^* \: + \:{\cal D} {\cal B}^* \: + \: 
{\cal B}^* {\cal A} \: + \: {\cal D}^* {\cal C} \notag \\ 
\delta \: & = \: {\cal D} {\cal D}^* \: + \: {\cal D}^* {\cal D} \: + \: 
{\cal B}^* {\cal B} \: + \: {\cal C} {\cal C}^* \: + \: k \: I. \notag
\end{align}

As $k \: > \: 0$, the operators $\alpha$ and $\delta$ are invertible,
with $\parallel \alpha^{-1} \parallel \:\le \: k^{-1}$ and
$\parallel \delta^{-1} \parallel \:\le \: k^{-1}$.  By Proposition
\ref{prop2}, $\delta^{-1}$ is compact. By an elementary argument, 
$\parallel {\cal D} \: \delta^{- \: 1/2} \parallel \: \le \: 1$,
$\parallel {\cal D}^* \: \delta^{- \: 1/2} \parallel \: \le \: 1$,
$\parallel \delta^{- \: 1/2} \: {\cal D} \parallel \: \le \: 1$ and
$\parallel \delta^{- \: 1/2} \: {\cal D}^* \parallel \: \le \: 1$.
Then ${\cal D} \: \delta^{- 1}$,
${\cal D}^* \: \delta^{-1}$,
$\delta^{- 1} \: {\cal D}$ and
$\delta^{- 1} \: {\cal D}^*$ are compact with norm at most $k^{- \: 1/2}$.

We claim that
$\alpha \: - \: \beta \: \delta^{-1} \: \gamma$ is invertible if $k$ is
large enough.
To see this, we write
\begin{equation}
\alpha \: - \: \beta \: \delta^{-1} \: \gamma \: = \: \alpha^{1/2} \: 
\left( I \: - \: \alpha^{- \: 1/2} \: 
\beta \: \delta^{-1} \: \gamma \: \alpha^{- \: 1/2} \right) \: \alpha^{1/2}.
\end{equation}
Then it suffices to show that 
$\parallel \alpha^{- \: 1/2} \: 
\beta \: \delta^{-1} \: \gamma \: \alpha^{- \: 1/2} \parallel \: < \: 1$ if
$k$ is large enough.
Writing things out, we have
\begin{align} \label{mess1}
\alpha^{- \: 1/2} \: 
& \beta \: \delta^{-1} \: \gamma \: \alpha^{- \: 1/2} \: = \\
& \alpha^{- \: 1/2} \: 
\left( 
{\cal A} {\cal C}^* \: + \: {\cal B} {\cal D}^* \: + \: 
{\cal A}^* {\cal B} \: + \: {\cal C}^* {\cal D} \right)
\: \delta^{-1} \: 
\left( {\cal C} {\cal A}^* \: + \:{\cal D} {\cal B}^* \: + \: 
{\cal B}^* {\cal A} \: + \: {\cal D}^* {\cal C} \right)
 \: \alpha^{- \: 1/2}. \notag
\end{align}
Now the operators $\alpha^{- \: 1/2} \: {\cal A}$, 
$\alpha^{- \: 1/2} \: {\cal A}^*$,
${\cal A}^* \: \alpha^{- \: 1/2}$, ${\cal A} \: \alpha^{- \: 1/2}$, 
${\cal D}^* \: \delta^{-1} \: {\cal D}$, ${\cal D} \: \delta^{-1} \: 
{\cal D}^*$, ${\cal D} \: \delta^{-1} \: {\cal D}$ and 
${\cal D}^* \: \delta^{-1} \: {\cal D}^*$ all have norm at
most one. From Proposition \ref{prop3},  ${\cal B}$ and ${\cal C}$ are bounded.
Writing out (\ref{mess1}) into its sixteen terms, we see that the
structure is such that by taking $k$ large, we can make the norm of
any individual term as small as desired. Hence by taking $k$ large
enough, we can make $\alpha \: - \: \beta \: \delta^{-1} \: \gamma$
invertible.

Writing 
\begin{equation} \label{mess2}
\left( \alpha \: - \: \beta \: \delta^{-1} \: \gamma \right)^{-1} \:
{\cal A} \: = \: \alpha^{- \: 1/2} \:  
\left( I \: - \: \alpha^{- \: 1/2} \: \beta \: \delta^{-1} \: \gamma \:
\alpha^{- \: 1/2} \right)^{-1} \:
\alpha^{- \: 1/2} \: {\cal A},
\end{equation}
we see that 
$\left( \alpha \: - \: \beta \: \delta^{-1} \: \gamma \right)^{-1} \:
{\cal A}$ is bounded. Similarly, 
$\left( \alpha \: - \: \beta \: \delta^{-1} \: \gamma \right)^{-1} \:
{\cal A}^*$, 
${\cal A} \: 
\left( \alpha \: - \: \beta \: \delta^{-1} \: \gamma \right)^{-1}$ and
${\cal A}^* \: 
\left( \alpha \: - \: \beta \: \delta^{-1} \: \gamma \right)^{-1}$
are bounded.
It now follows from (\ref{inverse}) that all components of
$\begin{pmatrix}
\alpha & \beta \\
\gamma & \delta
\end{pmatrix}^{-1}$ except for the upper left component
are compact. We note that the same statement is true
about $({\cal L} \: + \: k \: I)^{-1}$.
It remains to show that 
\begin{equation}
\left( \alpha \: - \: \beta \: \delta^{-1} \: \gamma
\right)^{-1} \: - \: \left( {\cal A} \: {\cal A}^* \: + \: {\cal A}^* \:
{\cal A} \: + \: k \: I \right)^{-1}
\end{equation}
 is compact.

Let us write 
\begin{equation} \label{mess3}
\alpha \: - \: \beta \: \delta^{-1} \: \gamma \: = \:
{\cal A} \: {\cal A}^* \: + \: {\cal A}^* \: {\cal A} \: + \: k \: I \:
- \: \left( \beta \: \delta^{-1} \: \gamma \: - \: {\cal B} \: {\cal B}^* \:
- \: {\cal C}^* \: {\cal C} \right).
\end{equation}
Then formally,
\begin{equation} \label{mess4}
\left( \alpha \: - \: \beta \: \delta^{-1} \: \gamma \right)^{-1} \: = \:
\left( {\cal A} \: {\cal A}^* \: + \: {\cal A}^* \: {\cal A} \: + \: k \: I 
\right)^{- \: 1/2} (I \: - \: X)^{-1} \: 
\left( {\cal A} \: {\cal A}^* \: + \: {\cal A}^* \: {\cal A} \: + \: k \: I 
\right)^{- \: 1/2},
\end{equation}
where
\begin{equation} \label{mess5}
X \: = \: 
\left( {\cal A} \: {\cal A}^* \: + \: {\cal A}^* \: {\cal A} \: + \: k \: I 
\right)^{- \: 1/2} \: 
\left( \beta \: \delta^{-1} \: \gamma \: - \: {\cal B} \: {\cal B}^* \:
- \: {\cal C}^* \: {\cal C} \right)
\left( {\cal A} \: {\cal A}^* \: + \: {\cal A}^* \: {\cal A} \: + \: k \: I 
\right)^{- \: 1/2}.
\end{equation}
It follows that 
\begin{align} \label{mess6}
& \left( \alpha \: - \: \beta \: \delta^{-1} \: \gamma \right)^{-1} \: - \:
\left( {\cal A} \: {\cal A}^* \: + \: {\cal A}^* \: {\cal A} \: + \: k \: I 
\right)^{- \: 1} \: = \\
&  \: \: \: \: \: 
\left( {\cal A} \: {\cal A}^* \: + \: {\cal A}^* \: {\cal A} \: + \: k \: I 
\right)^{- \: 1/2} \:
\left( \sum_{i=1}^\infty X^i \right)
\left( {\cal A} \: {\cal A}^* \: + \: {\cal A}^* \: {\cal A} \: + \: k \: I 
\right)^{- \: 1/2}, \notag
\end{align}
provided that the sum converges. We will show that $X$ is compact and that
the sum norm-converges if $k$ is large enough, which will prove the theorem.

We have
\begin{equation} \label{mess7}
\beta \: \delta^{-1} \: \gamma \: = \: 
\left( 
{\cal A} {\cal C}^* \: + \: {\cal B} {\cal D}^* \: + \: 
{\cal A}^* {\cal B} \: + \: {\cal C}^* {\cal D} \right)
\: \delta^{-1} \: 
\left( {\cal C} {\cal A}^* \: + \: {\cal D} {\cal B}^* \: + \: 
{\cal B}^* {\cal A} \: + \: {\cal D}^* {\cal C} \right).
\end{equation}
Consider first the terms of (\ref{mess7})
that are explicitly quadratic in ${\cal D}$, namely
\begin{equation} \label{mess8}
{\cal B} {\cal D}^* \: \delta^{-1} \: {\cal D} {\cal B}^* \: + \:
{\cal C}^* {\cal D} \: \delta^{-1} \: {\cal D}^* {\cal C} \: + \: 
{\cal B} {\cal D}^* \: \delta^{-1} \: {\cal D}^* {\cal C} \: + \:
{\cal C}^* {\cal D} \delta^{-1} \: {\cal D} {\cal B}^*.
\end{equation}
As $d^2 \: = \: 0$, we have ${\cal A} {\cal B} \: = \: - \: {\cal B}
{\cal D}$,
${\cal C} {\cal A} \: = \: - \: {\cal D} {\cal C}$ and
${\cal D}^2 \: = \: - \: {\cal C} {\cal B}$. 
Then (\ref{mess8}) equals
\begin{align} \label{mess9}
& {\cal B} {\cal D}^* \: \delta^{-1} \: {\cal D} {\cal B}^* \: + \:
{\cal B} {\cal D} \: \delta^{-1} \: {\cal D}^* {\cal B}^* \: + \:
{\cal C}^* {\cal D} \: \delta^{-1} \: {\cal D}^* {\cal C} \: + \: 
{\cal C}^* {\cal D}^* \: \delta^{-1} \: {\cal D} {\cal C} \: + \: \\
& {\cal B} {\cal D}^* \: \delta^{-1} \: {\cal D}^* {\cal C} \: + \:
{\cal C}^* {\cal D} \: \delta^{-1} \: {\cal D} {\cal B}^* \: - \:
{\cal A} {\cal B} \: \delta^{-1} \: {\cal B}^* {\cal A}^* \: - \:
{\cal A}^* {\cal C}^* \: \delta^{-1} \: {\cal C} {\cal A}. \notag
\end{align}
Thus
\begin{align} \label{mess10}
\beta \: & \delta^{-1} \: \gamma \: - \: {\cal B} \: {\cal B}^* \:
- \: {\cal C}^* \: {\cal C} \: = \\
& {\cal B} \left( {\cal D}^* \: \delta^{-1} \: {\cal D} \: + \:
{\cal D} \: \delta^{-1} \: {\cal D}^* \: - \: I \right) {\cal B}^* \: + \:
{\cal C}^* \left( {\cal D}^* \: \delta^{-1} \: {\cal D} \: + \: 
{\cal D} \: \delta^{-1} \: {\cal D}^* \: - \: I \right) {\cal C} \: + \: 
\notag \\
& {\cal B} {\cal D}^* \: \delta^{-1} \: {\cal D}^* {\cal C} \: + \:
{\cal C}^* {\cal D} \: \delta^{-1} \: {\cal D} {\cal B}^* \: + \:
O({\cal D}), \notag
\end{align}
where $O({\cal D})$ denotes the 
terms that are linear in ${\cal D}$.

We have
\begin{align} \label{mess11}
{\cal D}^* \: \delta^{-1} \: {\cal D} \: + \:
{\cal D} \: \delta^{-1} \: {\cal D}^* \: - \: I \: = \:
& \left( {\cal D}^* {\cal D} \: + \: {\cal D} {\cal D}^* \: - \: \delta \right)
\delta^{-1} \: + \\
& D^* \: \delta^{-1} ({\cal D} \: \delta \: - \: \delta \:
{\cal D} ) \: \delta^{-1} \: + \notag \\
& D \: \delta^{-1} ({\cal D}^* \: \delta \: - \: \delta \:
{\cal D}^* ) \: \delta^{-1}, \notag
\end{align}
\begin{align} \label{mess12}
{\cal D} \: \delta^{-1} \: {\cal D} \: & = \:
{\cal D}^2 \: \delta^{-1} \: + \: {\cal D} \: \delta^{-1} \: 
({\cal D} \: \delta \: - \: \delta \:
{\cal D} ) \: \delta^{-1} \\
& = \: - \: {\cal C} {\cal B} 
\: \delta^{-1} \: + \: {\cal D} \: \delta^{-1} \: 
({\cal D} \: \delta \: - \: \delta \:
{\cal D} ) \: \delta^{-1} \notag
\end{align}
and
\begin{align} \label{mess13}
{\cal D}^* \: \delta^{-1} \: {\cal D}^* \: & = \:
({\cal D}^*)^2 \: \delta^{-1} \: + \: {\cal D}^* \: \delta^{-1} \: 
({\cal D}^* \: \delta \: - \: \delta \:
{\cal D}^* ) \: \delta^{-1} \\
& = \: - \: {\cal B}^* {\cal C}^* 
\: \delta^{-1} \: + \: {\cal D}^* \: \delta^{-1} \: 
({\cal D}^* \: \delta \: - \: \delta \:
{\cal D}^* ) \: \delta^{-1}. \notag
\end{align}
Furthermore,
\begin{equation} \label{mess14}
{\cal D}^* {\cal D} \: + \: {\cal D} {\cal D}^* \: - \: \delta \: = \:
- \: {\cal B}^* \: {\cal B} \: - \: {\cal C} \: {\cal C}^* \: - \: k \: I,
\end{equation}
\begin{align} \label{mess15}
{\cal D} \: \delta \: - \: \delta \: {\cal D} \: & = \:
[ {\cal D}^2 , {\cal D}^*] \: + \: [{\cal D}, {\cal B}^* {\cal B} \: + \:
{\cal C} {\cal C}^*] \\
&= \: - \:  [ {\cal C} {\cal B} , 
{\cal D}^*] \: + \: [{\cal D}, {\cal B}^* {\cal B} \: + \:
{\cal C} {\cal C}^*] \notag
\end{align}
and
\begin{align} \label{mess16}
{\cal D}^* \: \delta \: - \: \delta \: {\cal D}^* \: & = \:
[ ({\cal D}^*)^2 , {\cal D}] \: + \: [{\cal D}^*, {\cal B}^* {\cal B} \: + \:
{\cal C} {\cal C}^*] \\
&= \: - \:  [ {\cal B}^* {\cal C}^* , 
{\cal D}] \: + \: [{\cal D}^*, {\cal B}^* {\cal B} \: + \:
{\cal C} {\cal C}^*]. \notag
\end{align}
Substituting (\ref{mess14}) - (\ref{mess16}) into 
(\ref{mess11}) - (\ref{mess13}), we see that
${\cal D}^* \: \delta^{-1} \: {\cal D} \: + \:
{\cal D} \: \delta^{-1} \: {\cal D}^* \: - \: I$,
${\cal D} \: \delta^{-1} \: {\cal D}$ and
${\cal D}^* \: \delta^{-1} \: {\cal D}^*$ are compact. 
Substituting (\ref{mess10}) into (\ref{mess5}), 
we see that the contributions to $X$ of the
terms listed in (\ref{mess10}) are all compact.

Next, from (\ref{mess7}), the terms in 
$\beta \: \delta^{-1} \: \gamma \: - \: {\cal B} \: {\cal B}^* \:
- \: {\cal C}^* \: {\cal C}$ that are explicitly proportionate to
${\cal D}$ are
\begin{align}
& {\cal A} {\cal C}^* \: \delta^{-1} \: {\cal D} {\cal B}^* \: + \: 
{\cal A} {\cal C}^* \: \delta^{-1} \: {\cal D}^* {\cal C} \: + \: 
{\cal B} {\cal D}^* \: \delta^{-1} \: {\cal C} {\cal A}^* \: + \: 
{\cal B} {\cal D}^* \: \delta^{-1} \: {\cal B}^* {\cal A} \: + \: \\ 
& {\cal A}^* {\cal B} \: \delta^{-1} \: {\cal D} {\cal B}^* \: + \: 
{\cal A}^* {\cal B} \: \delta^{-1} \: {\cal D}^* {\cal C} \: + \: 
{\cal C}^* {\cal D} \: \delta^{-1} \: {\cal C} {\cal A}^* \: + \: 
{\cal C}^* {\cal D} \: \delta^{-1} \: {\cal B}^* {\cal A}. \notag
\end{align}
One sees that their contributions to (\ref{mess5}) are all compact.
Finally, the remaining terms in 
$\beta \: \delta^{-1} \: \gamma \: - \: {\cal B} \: {\cal B}^* \:
- \: {\cal C}^* \: {\cal C}$ that are constant in ${\cal D}$ are
\begin{align}
\left( 
{\cal A} {\cal C}^* \: + \:
{\cal A}^* {\cal B} \right)
\: \delta^{-1} \: 
\left( {\cal C} {\cal A}^* \: + \: 
{\cal B}^* {\cal A} \right) \: - \:
{\cal A} \: {\cal B} \: \delta^{-1} \: {\cal B}^* {\cal A}^* \: - \: 
{\cal A}^* \: {\cal C}^* \: \delta^{-1} \: {\cal C} {\cal A}.
\end{align}
Their contributions to (\ref{mess5}) are all compact.

One can show as before
that the norm of $X$ can be made arbitrarily small by making $k$ large
enough.

This proves Theorem \ref{thm2} under our assumption that the metric on $M$ is
a product near each $\{0\} \times N_I$. However, as in
\cite[Proposition 2.1]{Donnelly-Li (1979)}, the essential
spectrum of $\triangle^M_p$ is invariant under a compactly-supported
change of the metric. Furthermore, 
the essential spectrum of a self-adjoint 
ordinary differential operator on $[0, \infty)$
is independent of the choice of (self-adjoint) boundary condition at
$\{0\}$ 
\cite[Volume II, Chapter XIII.7, Corollary 3]{Dunford-Schwartz (1963)}
and is also unchanged by a compactly-supported perturbation of the
operator.  Thus Theorem \ref{thm2} also holds for the original metric on $M$.

\section{Proof of Theorem \ref{thm3}} \label{sect8}

We now specialize to the case of functions.  In this case, 
$E_I$ is a trivial real line bundle on $[0, \infty)$. 
Consider the quadratic form (\ref{Q(F,F)}) in the case $B = 1$, with
$f \in C^\infty([0, \infty))$ and $f(0) \: = \: 0$.

Let $v(s)$ denote the volume of $(N, h(s))$. Then
\begin{equation}
\frac{dv}{ds} \: = \: - \: \int_{N} \sum_i S^i_{\: i} \: d\vol(s).
\end{equation}
If ${\cal F} : (0, \infty) \times N \rightarrow U$ is smooth then
the Gauss-Codazzi equation gives
\begin{equation}
- \: 
\partial_s \sum_i S^i_{\: i} \: + \: \sum_{ij} S^{ij} \: S_{ij} \: = \: - \: 
\Ric(\partial_s, \partial_s),
\end{equation}
which in turn implies that
\begin{equation}
\frac{d^2v}{ds^2} \: = \: \int_{N} \left[ - \: \Ric(\partial_s, \partial_s)
\: - \: \sum_{ij} S^{ij} S_{ij} \: + \: \left( \sum_i S^i_{\: i} \right)^2
\right]  \: d\vol(s).
\end{equation}
This last equation makes sense even if ${\cal F}$ is not smooth, showing
that $v$ is $C^2$-smooth in $s$.

\begin{lemma} \label{equiv}
$\left( {\cal A} \: {\cal A}^* \: + \: {\cal A}^* \: {\cal A} \right)_0$ 
is unitarily equivalent to the operator
\begin{equation}
- \: \frac{d^2}{ds^2} \: + \: 
\frac{1}{2} \: \frac{d^2\ln v}{ds^2} \: + \: \frac{1}{4}
\left( \frac{d\ln v}{ds} \right)^2,
\end{equation}
which is densely-defined and self-adjoint on $L^2([0, \infty))$, with
Dirichlet boundary conditions.
\end{lemma}
\begin{pf}
Putting
$k(s) = v(s)^{1/2} f$, we have
\begin{equation}
\langle f, f \rangle = \langle k, k \rangle_{L^2}
\end{equation}
and
\begin{align}
Q(f) & = \int_{0}^\infty 
\left( \frac{d}{ds} ( v^{-1/2} k) \right)^2 \: v(s)
\: ds \\
& = \int_{0}^\infty \left( v^{-1/2} \frac{dk}{ds} \: - \: \frac{1}{2} \:
v^{-3/2} \: \frac{dv}{ds} \: k \right)^2 \: v(s) \: ds \notag \\
& = \int_{0}^\infty \left( \frac{dk}{ds} \: - \: \frac{1}{2} \:
v^{-1} \: \frac{dv}{ds} \: k \right)^2 \: ds \notag \\ 
& = \int_{0}^\infty \left[ \left( \frac{dk}{ds} \right)^2 \: - \:
v^{-1} \: \frac{dv}{ds} \: k \: \frac{dk}{ds} \: + \: \frac{1}{4}
\left( v^{-1} \: \frac{dv}{ds} \right)^2 \: k^2
\right] \: ds \notag \\  
& = \int_{0}^\infty 
\left[ \left( \frac{dk}{ds} \right)^2 \: - \: \frac{1}{2} \:
\frac{d\ln v}{ds} \: \frac{dk^2}{ds} \: + \: \frac{1}{4}
\left( \frac{d\ln v}{ds} \right)^2 \: k^2
\right] \: ds \notag \\  
& = \int_{0}^\infty \left[ \left( \frac{dk}{ds} \right)^2 \: + \: 
\left( \frac{1}{2} \: \frac{d^2\ln v}{ds^2} \: + \: \frac{1}{4}
\left( \frac{d\ln v}{ds} \right)^2 \right) \: k^2
\right] \: ds. \notag 
\end{align}
The lemma follows.
\end{pf}

Now let $P$ be an even periodic element of $C^\infty(\R)$ which is
not real-analytic. Put 
\begin{equation}
V_P \: = \: \frac{1}{2} \: \frac{dP}{ds} \: + \: \frac{1}{4} \:
P^2. 
\end{equation}
Let ${\cal O}_P^\prime$ 
be the operator $- \: \frac{d^2}{ds^2} \: + \: V_P$ acting on
$L^2([0, \infty))$, with Dirichlet boundary conditions at $0$.

\begin{lemma} \label{gaps}
${\cal O}_P^{\prime}$ has an infinite number of gaps in its
essential spectrum.
\end{lemma}
\begin{pf}
As $\frac{dP}{ds}$ is odd and $P^2$ is even, if $V_P$ were real-analytic
then $\frac{dP}{ds}$ would be real-analytic, which would imply that
$P$ is real-analytic.  Thus $V_P$ is not real-analytic. From
\cite[Vol. IV, Thm. XIII.91(d)]{Reed-Simon (1978)}, the operator 
${\cal O}_P \: = \: - \: \frac{d^2}{ds^2} \: + \: V_P$
on $L^2(\R)$
has an absolutely continuous spectrum which consists of an infinite
number of disjoint closed intervals in $[0, \infty)$, tending toward
infinity. Let
${\cal O}_P^{\prime \prime}$ be the operator 
$- \: \frac{d^2}{ds^2} \: + \: V_P$ acting on
$L^2((-\infty, 0])$, again with Dirichlet boundary conditions at $0$.
Then the essential spectrum of ${\cal O}_P$ is the union of the essential
spectra of ${\cal O}_P^{\prime}$ and ${\cal O}_P^{\prime \prime}$. 
As the essential spectra of both
${\cal O}_P^{\prime}$ and ${\cal O}_P^{\prime \prime}$ tend toward
infinity, the lemma follows.
\end{pf}
\noindent
{\bf Proof of Theorem \ref{thm3} :} \\

Start with a complete finite-volume hyperbolic metric on a punctured
$2$-torus.  On the cusp, the metric is $ds^2 + e^{-2s} d\theta^2$ for
$s \in [s_0, \infty)$, with $s_0 \:  > \: 0$. 

Let $p$ be an even periodic element of $C^\infty(\R)$ which is not
real-analytic. Let $\phi \in C^\infty_0([0, \infty))$ be a nonincreasing
function which is identically one on $[0,1]$ and identically zero on
$[2, \infty)$. For $\delta \: > \> 0$ and $s \: \ge \: s_0$, put 
\begin{equation}
v_\delta(s) 
\: = \: e^{-s - \delta \int_0^{s-s_0} p(u) \: (1 - \phi(\delta u)) \: du}.
\end{equation}
Keep the metric on the complement of the cusp unaltered and
change the metric on the cusp to $ds^2 \: + \: v_\delta(s)^2 \:
d\theta^2$. From Theorem \ref{thm2}
and Lemma \ref{equiv}, the essential spectrum of the
Laplacian of the new metric is the same as the essential spectrum of the 
operator ${\cal O}_{-1 - \delta p}^\prime$. Then from Lemma
\ref{gaps}, the Laplacian of the new metric
has an infinite number of gaps in its essential spectrum.
Hence it has an infinite number of gaps in its spectrum.  One
can check that as $\delta \rightarrow 0$, the sectional curvatures of the new
metric become pinched arbitrarily close to $-1$. \\ \\
{\bf Remark : } It seems likely that by taking $p$ to be almost-periodic,
one can find similar examples in which the essential spectrum is a Cantor set.

\end{document}